\documentclass[12pt]{amsart}
\usepackage{amsmath,amssymb,amsthm,amscd,amsopn,amsfonts,amsxtra}
\usepackage{latexsym,array}
\usepackage[mathscr]{eucal}
\usepackage{graphicx,psfrag}
\usepackage{epsfig}
 \newtheorem{theorem}[subsection]{Theorem}
\newtheorem{lemma}[subsection]{Lemma}

\newtheorem{proposition}[subsection]{Proposition}
\newtheorem{definition}[subsection]{Definition}
\newtheorem{remark}[subsection]{Remark}

\def\al{\alpha}

\def\eps{\epsilon}

\def\vphi{\varphi}
\def\lam{\lambda}

\def\sh{\sinh}
\def\ch{\cosh}
\def\th{\tanh}
\def\sech{\mathrm{sech}}
\def\supp{\mathrm{supp}}

\newcommand{\cal}{\mathcal}

\newcommand{\nd}{\noindent}

\newcommand{\n}{\newline}

\newcommand{\Z}{\mathbb{Z}}
\newcommand{\R}{\mathbb{R}}
\newcommand{\C}{\mathbb{C}}
\newcommand{\N}{\mathbb{N}}

\numberwithin{equation}{section}

\begin{document}
\title[Function space associated with Schr\"odinger operators]{Function spaces
associated with Schr\"odinger operators: the P\"oschl-Teller potential}
\author{Gestur \'Olafsson}
\address[Gestur \'Olafsson]{Department of Mathematics \\
         Louisiana State University  \\
         Baton Rouge, LA 70803}
 \email{olafsson@math.lsu.edu}
           \urladdr{http://www.math.lsu.edu/\textasciitilde{olafsson}}
\author{Shijun Zheng}
\address[Shijun Zheng]{Department of Mathematics \\
                       Industrial Mathematics Institute\\
                     University of South Carolina\\    
                     Columbia, SC 29208}
\address{  and}
\address{Department of Mathematics\\
Louisiana State University  \\
         Baton Rouge, LA 70803}
\email{shijun@math.sc.edu}
 \urladdr{http://www.math.sc.edu/\symbol{126}{shijun}}
\thanks{G. \'Olafsson was supported by NSF grants
DMS-0139473 and DMS-0402068.
S. Zheng was partially supported by
DARPA grant MDA 972-01-1-0033. He also gratefully thanks
NSF for support during his visit at Louisiana State University in 2003}

\date{\today}

\begin{abstract}
We address the function space theory associated with the Schr\"odinger operator
 $H=-d^2/dx^2+V$.  The discussion is featured with 
potential
 $V(x)= -\lambda(\lambda-1)\,\sech^2 x $, 
which is called in quantum physics the P\"oschl-Teller potential. 
Using biorthogonal  dyadic system, we introduce Besov spaces and
Triebel-Lizorkin spaces (including Sobolev spaces) associated with $H$. 
We then use interpolation method to identify these spaces  
with the classical ones for a certain range of $p,q> 1$.  
A physical implication is that the corresponding 
wave function $\psi(t,x)=e^{-i t H }f(x)$ admits appropriate 
time decay 
in the Besov space scale. 
\end{abstract}
 
\maketitle

\section*{Introduction}
\noindent
Let $H = - d^2/dx^2 +V$ be a
Schr\"{o}dinger operator on $\R$ with  real valued potential function
$V$. 
In quantum physics, $H$ is the energy
operator of a particle with
one degree of freedom with potential $V(x)$. The Schr\"odinger equation 
$$ i \; \partial_t \psi= H \psi,\qquad
\psi(0,x)=f(x)$$
has solution given by  $\psi(t,x)=e^{-i t H }f(x)$.
If the potential has
certain decay at $\infty$, then one expects that asymptotically, as time
tends to infinity, the motion of the associated perturbed quantum system
resembles the free evolution.
In fact, if $V$ is in $L^1 (\R) \cap L^2 (\R )$
or if $\int_\R (1+|x|^2)|V(x)|dx  < \infty$,
then it is known, c.f.  \cite{CK01}, [GH98], \cite{Z04a}, 
that the essential spectrum of $H$ and $H_0=-d^2/dx^2 $ coincide and there is no
singular continuous spectrum; the wave operators
$W_{\pm}=s-\lim_{t\rightarrow \pm \infty}e^{ i t H } e^{-i t H_0}$ exists and are
complete.

Recently, several authors have considered
function spaces associated with Schr\"odinger operators, cf. \cite{JN94, E95,
E96, DZ98, DZ02, BZ04}.
One of the goals has been to develop the associated Littlewood-Paley theory,
in order to give a unified approach.
Motivated by the treatment in \cite{BZ04,E95}
we study the negative
potential
\begin{equation}\label{eq-potential}
V(x)= -\lambda(\lambda-1) \sech^2 x\,
\end{equation}
called P\"oschl-Teller potential. It
is a fundamental model in soliton theory.
In this article, we
will mainly consider the case where $\lam=n+1$, $n\in \N_0$. In that case
we obtain a compact expression for eigenfunctions of $H$.
It turns out that for the absolute continuous part of $H$,
the analysis is simpler than the barrier potential,
although $H$ has a nonempty pure point spectrum.

Denote by $C_0(\R )$ the space of compactly supported functions
on the line.
Suppose $\Phi, \varphi \;\Psi, \psi\in C_0^\infty ({\R}) $ satisfy
the ``biorthogonal'' condition (see Section \ref{S3})
\begin{equation}\label{eq:dyadic-id}
\Phi(x) \Psi(x) +\sum_{j=1}^\infty \varphi_j(x) \psi_j(x) =1,
\end{equation}
where $\varphi_j(x) =\varphi( 2^{-j} x)$.
Let $\alpha \in \mathbb{R}, 0<p,q< \infty$. The {\em Triebel-Lizorkin spaces}
associated with $H$, denoted by
$F_p^{\alpha,q}:=F_p^{\alpha,q}(H)$,
is defined to be the completion of the subspace $ L_0^2(\R ):= \{f\in L^2(\R ):
\Vert f\Vert_{F_p^{\alpha,q}} <\infty \}$,
where the quasi-norm $\Vert \cdot\Vert_{F_p^{\alpha,q}}$ is initially
defined for $f\in L^2(\R )$ by:
\begin{equation}\label{eq:F-norm}
\Vert f\Vert_{F_p^{\alpha,q}} 
=\Vert \bigg(\sum_{j=0}^{\infty} 2^{j\alpha q} \vert \vphi_j(H)f \vert^q\bigg)^{1/q}\Vert_p.
\end{equation}

The main result is an equivalence theorem for $F_p^{\alpha,q}(H)$
using Peetre's maximal function, which suggests that
that $\Vert f\Vert_{F_p^{\alpha,q}}^{\phi}
$ and $\Vert f\Vert_{F_p^{\alpha,q}}^{\psi}$ are equivalent
quasi-norms on $F_p^{\alpha,q}(H)$, for given
two systems of smooth dyadic systems $\{\phi_j\}, \{\psi_j\}$.
The following is proved in Theorem \ref{th:phi*F-inhomo}:

\begin{theorem}  Let $\alpha\in \R,  0<p,q \leq \infty$.
If $\vphi_j^*f$
is defined for $j\geq 0 $ with $s>1/\min(p,q)$, we have for $f\in L^2 (\R)$
\begin{equation}
 \Vert f\Vert_{F_p^{\alpha,q}} \sim \Vert {\Phi}^*f \Vert_p+
\Vert
\left(\sum_{j=1}^{\infty} ( 2^{j\alpha} \vphi_j^*f)^q \right)^{1/q} \Vert_p.
\end{equation}
Furthermore, $F_p^{\alpha,q}$ is a quasi-Banach space (Banach space if
$p\geq 1, q \geq 1$) and it is independent of the choice of $\{\Phi,
{\vphi}_j\}_{j\geq 1}$.
\end{theorem}

Here $\vphi_j^*f=\vphi_{j,s}^*f$ denotes the Peetre type maximal function, see
Section \ref{S3} for details.


As in \cite{BZ04}
the essential part is to obtain the decay estimates of the kernel of
the operator
$\phi_{\lambda}(H)$ as well as of its derivatives.  This requires an
examination of the differential properties of the eigenfunctions
$ e_k(x):=e(x, k)$. For $k=i, 2i, \dots, ni$, where
as usually $i=\sqrt{-1}$, the bounded states are Schwartz functions.
In Section \ref{S2} we solve the Lippmann-Schwinger equation (\ref{e-ls})
for $k\in \R\cup \{i,\dots, ni\}$ based on a general expression
suggested in
\cite[p.40, p.54]{Lam80}. Furthermore, we determine precisely the transmission
and reflection coefficients.

We mention that our results are also true for translation and regular
scaling (i.e.,
$V\longrightarrow \alpha^2 V(\alpha)$) of the potential $V$.

A natural question is: what is the relation between the new function spaces
and classical ones, namely, $F^{\alpha,q}_p(\R)$ and $B^{\alpha,q}_p(\R)$.
Concerning identification of the $F$ and $B$ spaces, it
can be shown that the following holds:

\nd
{\bf Theorem \ref{th:F-interp}\;\;} If $\al>0$, $1<p<\infty$
and $2p/(p+1)<q<2p$, then
$$F^{\alpha, q}_p(H) = F^{2\alpha,q}_p(\R)$$
and if $\al>0$, $1\le p<\infty$, $1\le q\le \infty$, then
$$B^{\alpha, q}_p(H) = B^{2\alpha,q}_p(\R).$$

For $0<p\le 1$, it would be interesting to see whether it holds that
$F^{0,1}_p={\cal H}^p$, where ${\cal H}^p$ is the Hardy spaces associated with $H$
defined
by atomic decomposition or the maximal function associated to the
heat kernel, cf. \cite{DZ98, DZ02}.

In a sequal to this paper we will study the
spectral multiplier problem on the $F$ and $B$ spaces.

\section{The eigenfunctions of H}\label{S2}

In this section we derive a simple expression for
the eigenfunctions on the Schr\"odinger operator
$H=-d^2/dx^2-\lambda (\lambda -1)\sech^2(x)$
in the case wherer $\lambda = n+1$ is an integer.

\subsection{The Lippmann-Schwinger equation}
Consider the eigenvalue problem
\begin{equation}\label{e-evp}
He(x,k)=k^2e(x,k),
\end{equation}
with asymptotic behavior
\begin{equation}\label{e-epm}
e_\pm (x,k)\sim \left\{\begin{array}{c@{\quad \mathrm{for}\quad}l}
T_{\pm}(k)e^{ikx} & x\to \pm \infty\\
e^{ikx} + R_\pm (k)e^{-ikx} & x\to \mp \infty\end{array}
\right.
\end{equation}
We will often use the notation
\begin{equation}\label{e-e}
e(x,k)= \left\{\begin{array}{c@{\quad \mathrm{for}\quad}l}
e_+(x,k) & k>0\\
e_-(x,k) & k<0\end{array}
\right.
\end{equation}
The coefficients $T_\pm (k)$ and $R_\pm (k)$
in (\ref{e-epm} are called the \textit{transmission
coefficients} and \textit{reflection coefficients}, resp.
They
satisfy the relation
\begin{equation}\label{e-coef}
|T_\pm(k)|^2 + |R_\pm (k)|^2=1\, .
\end{equation}
The equations (\ref{e-evp}) and (\ref{e-epm}) are
equivalent to the Lippmann-Schwinger equation
\begin{equation}\label{e-ls}
e_\pm (x,k)=e^{ikx}+\frac{1}{2i|k|}\int e^{i|k|\, |x-y|}V(y)e_\pm
(y,k)\, dy,
\end{equation}
where $\pm $ indicates the sign of $k$.
For
$ V(x)= -\lambda(\lambda-1) \sech^2 x, $ we seek a solution to
the integral equation (\ref{e-ls}) for each $k$.

\subsection{Inductive construction of the solution}\label{ss-ics}
Let us recall how to obtain the solution
for the eigenvalue problem (\ref{e-evp})
for integer values of $\lambda=n+1$.
By regularity, a weak solution $u$ solving $Hu=k^2 u$
is smooth on the domain where $u$ is smooth.
For $n\in \mathbb{N}_0$ denote by $T_n$
the differential operator
\begin{equation}\label{d-Tn}
T_n=\frac{d\, }{dx}-n \tanh x\, .
\end{equation}
With $y_n=e(x,k)$ the equation (\ref{e-evp}) can then be
written as
\begin{equation}\label{e-Dn}
y_n''+  n(n+1)\sech^2 x\, y_n= -k^2 y_n\, .
\end{equation}
If $n=0$, then the general solution is given
by $y_0=A e^{ik x}+B e^{-ik x}$. For $n \ge 1$ we obtain
\begin{equation}\label{e-induction}
y_{n}(x) = T_ny_{n-1}(x) = T_n\cdots T_1y_0(x)
=A(k)e^{ikx}+B(k)e^{-ikx}\, .
\end{equation}
\begin{lemma}\label{le-AB}
Assume that $\lambda =n+1$, with $n\in \N_0$.
Then there exists a polynomial $p_n(t,ik)$
of degree $n$ in both $t$ and $k$, such that
$$e_\pm (x,k)=A_n^\pm(k)p_n^\pm (\th x,ik)e^{ikx}\, .$$
Furthermore the following holds:
\begin{enumerate}
\item The constant $A_n^\pm (k)$ is given by
$$ A_n^+(k)=\prod_{j=1}^n\frac{1}{j+ ik}\qquad\mathrm{and}
\qquad A_n^-(k)=(-1)^n\prod_{j=1}^n\frac{1}{j-ik}\, .$$
\item The transmission coefficient $T_\pm (k)$
are
$$T_+(k)=(-1)^{n}\prod_{j=1}^n\frac{j-  ik}{j+  ik}
\qquad\mathrm{and}
\qquad T_-(k)=(-1)^n\prod_{j=1}^n\frac{j+ik}{j-ik}\, .$$
\item The reflection coefficients $R_\pm (k)$ are
all zero.
\end{enumerate}
\end{lemma}
\begin{proof}
As $n$ is fixed here, we will omit the subscript $n$ in the
following except in obvious induction steps.
As $\frac{d}{dx}(\th x) =1-\th^2 x$ it follows  from (\ref{e-induction})
that
\begin{equation}\label{eq:y_n}
e_\pm (x,k)=A^\pm (k)p^\pm_n ( x,ik)e^{ikx}+B(k)q^\pm_n( x,ik)e^{-ikx},
\end{equation}
where $p^\pm (x,k):=p^\pm(\th x, ik)$
and $q^\pm (x,k)=q^\pm(\th x,ik)$
are polynomials of degree $n$ in $k$ and $\th x$.

Assume for the moment, that $k>0$. We will
therefore also omit the superscript $\pm$
for the moment.
Substituting (\ref{eq:y_n}) into
 the asymptotics for $k>0$,\\
$$e(x,k)\longrightarrow \left\{
\begin{array}{ll}
 T_+(k)e^{ikx}, &x\rightarrow \infty,\\
 e^{ikx}+ R_+(k)e^{-ikx}, &x\rightarrow -\infty.
\end{array}\right.
$$
we obtain  $B(k)=0=R_+(k)$,
\begin{equation}\label{eq-Ak}
A(k)p(-1,ik)=1
\end{equation}
and
\begin{equation}\label{eq-Tk}
T_+(k)=A(k)p(1,ik)\, .
\end{equation}
This implies easily the recurrence formula
\begin{equation}\label{eq:p-recur}
p_n (\th x, ik)=
\frac{d\, }{dx}p_{n-1}(\th x, ik)+( ik-n\th x) p_{n-1}(\th x, ik).
\end{equation}
Note that $\frac{d\, }{dx}p_{n-1}(\th x, ik)=
p_{n-1}^\prime(\th x,ik)(1-\th^2(x))$. As
$p_{n-1}^\prime (t,ik)$ is a polynomial in $t$ it follows that
$$\lim_{x\to \pm \infty}p_{n-1}^\prime(\th x,ik)=
p_{n-1}^\prime (\pm 1,ik)$$
is bounded, and hence
\begin{equation}\label{eq-0}
\lim_{x\to \pm \infty }\frac{d\, }{dx}p_{n-1}(\th x,ik) =0\, .
\end{equation}
By the recursion formula (\ref{eq:p-recur}) we therefore get
\begin{equation}\label{eq:pn(-1)}
p_n(\pm 1,ik)=( ik\mp n)p_{n-1}(\pm 1,ik)
\end{equation}
and by induction
\begin{equation}\label{eq-p1}
p(1,ik)=(-1)^n\prod_{j=1}^n j-ik
\end{equation}
and
\begin{equation}\label{eq-p(-1)}
p(-1,ik)= \prod_{j=1}^n j+ik=(-1)^n\overline{p(1,k)}
\end{equation}
By (\ref{eq-Ak}) we now get:
\begin{equation}\label{eq-Ak1}
A_n(k)=\frac{1}{p_n(-1,ik)}=\prod_{j=1}^n\frac{1}{j+ik}\, .
\end{equation}

By the recursion formula (\ref{eq:p-recur}) and (\ref{eq:pn(-1)})
we get $p(1,ik)=\prod_{j=1}^n(j-ik)^{-1}$ and hence
(\ref{eq-Ak}) and
(\ref{eq-Tk}) imply the the statement for $T_+(k)$.

For $k$ negative we notice the following changes have
to be made. Instead of (\ref{eq-Ak}) we now have
\begin{equation}\label{eq-Ak-}
A_n(k)p(1,ik)=1
\end{equation}
and instead of (\ref{eq-Tk}) we have
\begin{equation}\label{eq-Tk}
T_+(k)=A(k)p(-1,ik)\, .
\end{equation}
The claim follows now from (\ref{eq-p1}) and
(\ref{eq-p(-1)}).
\end{proof}

The proof shows that $p_n^+(\th x,ik)=p_n^-(\th x,ik)$. Furthermore, simple
induction shows that
\begin{equation}\label{eq-psignchagen}
p_n(\th x,-k)=(-1)^n p(\th (-x),ik)\, .
\end{equation}
We therefore get
the following explicit formula for the eigenfunctions $e(x,k)$:
\begin{theorem}\label{th-exk}
Assume that $k\not= 0$. Then
$$e(x,k)=(\mathrm{sign}(k))^{n}\left(\prod_{j=1}^n\frac{1}{j+i | k|}
\right)\,P_n(x,ik)e^{ikx}$$
where $P_n(x,k)=p_n(\th (x),ik)$ is defined by the recursion formula
$$
p_n (\th x, ik)=
\frac{d\, }{dx}p_{n-1}(\th x, ik)+( ik-n\th x) p_{n-1}(\th x, ik)\, .$$
In particular the function
$$\R\times (\R\setminus \{0\})\ni (x,k)\mapsto e(x,k)\in\C$$
is analytic, and $e(x,-k)=e(-x,k)$. Finally, the function
$$(x,y,k)\mapsto e(x,k)\overline{e(y,k)}
=\left(\prod_{j=1}^n\frac{1}{j^2+ k^2}
\right)\, P_n(x,ik)P_n(-y,ik) e^{ik(x-y)}$$
is real analytic on $\R^3$.
\end{theorem}

\subsection{The point spectrum} In this subsection
we discuss briefly the point spectrum
of $H$. We start with:

\begin{lemma} If $\lambda=n(n+1)$, then the point
spectrum of $H$ is given by
$$\sigma_p=\{-1,-4,\ldots ,-n^2\}\, .$$
\end{lemma}
\begin{proof}
By \cite[Theorem 6.1, Theorem 7.3]{GH98} or \cite{Ti58}
the point spectrum is given by
the points $\mu^2$, such that
$\mathrm{Im}\mu > 0$ and $\mu$ and $T(k)$ has a pole
at $\mu$. The claim now follows from the explicite
description of $T(k)$ in
Lemma \ref{le-AB}.
\end{proof}

In (\ref{e-Dn}) replace $k$ by $ji$, where
$j\in \{1,2,\dots,n\}$, and denote the solution
by $y_{n,j}$. We have two particular
solutions: $y_{n,n}= \cosh^{n} x$  and  $y_{n,n+1} = \sech^n x$.
Induction shows that
$$ y_{n,n+2}= T_{n+1} \sech^n x$$
and
$$ y_{n,n+3}= T_{n+2}y_{n,n+2}\, .$$
For $n=1$ the particular solution is $y_{1,3}= -3 \th x \,\sech x$.
For $n=2$ the point spectrum is $\sigma_p =\{-1, -4 \}$,
and the corresponding
eigenfunctions (bound states)
are $\sech x \th x$ and $\sech^2 x$.

Finally we note that (\cite[p.46]{Lam80})
 the bound-state solution is orthogonal to
the propagating solutions, that is,
$$ \int \sech (k z)  (i\eta -k \th (k z ) )e^{i\eta z}\, dz=0$$



\subsection{Examples}
Assume first that $\lambda=n+1=2$.
In this case $n=1$, the continuous spectrum is
$\sigma_c=[0,\infty)$, and
the point spectrum is $\sigma_p={-1}$.

The boundary condition $u(\pm\infty)=0$ leads to the
single eigenvalue $E= -1$ with the associated eigenfunction
$u= \sech (x)$.
In quantum theory, the interpretation of this result is that  a
particle is confined by a potential well having a shape proportional
to $\sech^2( x)$ while the single value of $E$ is proportional to the
energy that the particle confined by this well can possess.
When the sign of the potential is reversed so that $V(x)= 2 \sech^2 (x)$,
the potential is repulsive and no bound state occurs.
The solutions for the eigenfunctions in the continuous spectrum
are given by
$$e(x,k)= \mathrm{sign} (k)\frac{ik-\th (x)}{1+i|k|}\, e^{ikx}\, .$$

If $\lam =3$. Thus $n=2$. Then $\sigma_c=[0,\infty)$,
and $\sigma_p=\{-1, -4\}$.
The eigenfunctions in the continuous spectrum for$k\not= 0$
are give by

$$e (x,k)=\frac{3 \th^2(x)-3ik\th (x)- k^2-1}{- k^2+3i|k|+2}\, e^{-kx}\, .$$

For point spectrum $\sigma_p =\{-1, -4 \}$, the corresponding
eigenfunctions (bounded states)
are $\sech x \th x$ and $\sech^2 x$ respectively.

\subsection{The continuous parameter $\lam >1$}
For $\alpha >0$ set
$a=\frac{1}{2}(\lam +i
\frac{k}{\al} )$, $b =\frac{1}{2}(\lam -i
\frac{k }{\al} ) $.  Then the fundamental system
(\ref{e-evp}) of two
real standard solutions $u_e$ and $u_o$, 
even and odd in $x$, are given
by \cite{Flu74},
$$
u_e(x)= \ch^\lam (\al x) F(a,b,\frac{1}{2}; -\sh^2 \al x)$$
and
$$ u_e(x)= \ch^\lam (\al x)  \sh (\al x) F(a+\frac{1}{2},
b+\frac{1}{2},\frac{3}{2}; -\sh^2 \al x).$$

Since $-\sh^2 \al x \sim -\frac{e^{2\al |x| } }{4}$
and $|x|\to \infty$, we obtain by the
asymptotic of hypergeometric functions:
$$ u_e(x)\sim  \frac{e^{\lam \al |x| }\Gamma (\frac{1}{2})}{2^{\lambda}}
\left\{ \frac{\Gamma(b-a) }{ \Gamma(b ) \Gamma(\frac{1}{2}-a)} 2^{2a}e^{-2a
\al |x| }
+ 
\frac{\Gamma(a-b) }{ \Gamma(a ) \Gamma(\frac{1}{2}-b) }2^{2b}e^{-2b
\al |x|} \right\}
$$
and
\begin{eqnarray*}
 u_o (x) &\sim&
 \pm  \frac{e^{(\lam +1)\al |x| }\Gamma (\frac{3}{2})}{2^{\lam+1}}
\left\{ \frac{\Gamma(b-a) }{ \Gamma(b+\frac{1}{2} ) \Gamma(1-a)}
2^{2a+1}e^{(-2a+1)
\al |x| }\right.\\
& &\, \,  + \left.
\frac{\Gamma(a-b) }{ \Gamma(a +\frac{1}{2}) \Gamma(1-b) }
2^{2b+1}e^{-(2b+1) \al |x|}\right\}
\end{eqnarray*}
($\pm$ sign corresponding to $x>1 $ or  $x<1$), \cite[p.249]{Le75}.


If the energy $E=k^2$ is non-zero, which corresponds to the scattering solution,
we have the following asymptotic expressions for $x>0$ or $x<0$:
$$u_e(x)\simeq C_e \cos (k|x|+\phi_e )$$
and
$$u_o(x)\simeq \pm C_o \cos (k|x|+\phi_o)$$
with the abbreviations,
\begin{eqnarray}
\phi_e(x) &=&\arg \frac{\Gamma(ik/\al)e^{-i\frac{k}{\al}\log 2 } }{
\Gamma(\frac{\lam}{2}+ i\frac{k}{2\al} ) \Gamma(\frac{1-\lam }{2} +
i\frac{k}{2\al} )}\label{eq:phase-e}\\
\phi_o (x) &=&\arg \frac{\Gamma(ik/\al) e^{-i\frac{k}{\al}\log 2 } }{
\Gamma(\frac{\lam-1}{2}+ i\frac{k}{2\al} ) \Gamma(1-\frac{\lam }{2} +
i\frac{k}{2\al} )}\label{eq:phase-o}.
\end{eqnarray}

Seeking a solution of the asymptotic form
$$
u=\left\{
\begin{array}{ll}
Te^{ikx}& x>0\\
e^{ikx}+Re^{-ikx}& x<0
\end{array}\right.
$$
we compute $T, R$ as
$$T=(e^{2i\phi_e}-e^{2i\phi_o})/2,
R=(e^{2i\phi_e}+e^{2i\phi_o})/2.
$$
Hence,
$$
|T|^2=\sin^2(\phi_e-\phi_o), \;|R|^2=\cos^2(\phi_e-\phi_o),
$$
which satisfy the conservation law  $|T|^2+|R|^2=1$.

Using (\ref{eq:phase-e}),(\ref{eq:phase-o}) we  can write
$$
\phi_e-\phi_o=\tan^{-1}\left(\frac{\sinh \pi k}{\sin \pi \lam}\right).
$$

For the discrete spectrum let $k=ij$, $j\in \N_0$ yields the
eigenvalues $-(\lam-j-1)^2$, $j\le \lam -1$, $j\in \N_0$.
Hence for the eigenstate equation
$$-\frac{\hbar^2}{2m}\frac{d^2u}{dx^2}+V=Eu,$$
the corresponding energies are, for $j\le \lam-1$
$$E_n=-\frac{\hbar^2}{2m}(\lam-1-j)^2,
$$
where the zero energy is a resonance (as $k\to 0$) rather than
an eigenvalue.

\subsection{Projection of the spectral operator  $\phi(H)$}
Let ${\cal H}=L^2$.
Then ${\cal H}={\cal H}_{ac}\oplus {\cal H}_{pp}$,
where ${\cal H}_{ac}$ denotes the absolute continuous subspace of $H$ and
${\cal H}_{pp}$ the pure point subspace.
The corresponding orthogonal projections
will be denoted by $E_{ac}$ respectively $E_{pp}$. Note
that $E_{ac}=E_{[0,\infty)}$ and $E_{pp}=E_{\sigma_p}$,  where
$E$ stand for the spectral projection.
Hence, if
$f\in L^2$, then $f=E_{ac}f + E_{pp} f$. It follows that
$$
 \phi(H)f =
\phi(H) E_{ac} f + \phi(H)E_{pp} f
=\phi(H)\big\vert_{{\cal H}_{ac}} f + \phi(H)\big\vert_{{\cal H}_{pp}} f.
$$
For $\phi\in C_0$ we have
\begin{equation}\label{eq:phi(H)f}
\phi(H)f(x)=
\int K(x,y) f(y) dy +  \sum_{k^2 \in \sigma_p } \phi(k^2)(f, e_k) e_k
\end{equation}
where
\begin{equation}\label{eq:kernel-repre}
K(x,y)= (2\pi)^{-1} \int \phi(k^2)e(x,k)\bar{e}(y,k) dk
\end{equation}
 is the kernel of $\phi(H)\vert_{{\cal H}_{ac}}$.
As $e(\cdot ,k)$ is smooth for $k\not= 0$, it follows
from (\ref{eq:kernel-repre}) that $K(x,y)$ is smooth if $\phi$ is
compactly supported.
Hence, if $\sigma_p$ is finite, which is the case when $\int x|V(x)|dx<\infty$,
and if $e(\cdot, k)\in L^2 (\R) \cap L^p(\R )$, as in our case, we  have
$\Vert \phi(H)\big\vert_{{\cal H}_p} f \Vert_{L^p} \le C \Vert f
\Vert_{L^p}$,  $f \in L^2 \cap L^p$.

\section{Decay estimates of the spectral operators $\phi_j(H)$ }
In this section we introduce the Triebel-Lizorkin and Besov spaces
associated with the Schr\"odinger operator
$H$ and show that they are well-defined,
i.e.,  different dyadic systems give rise to equivalent quasi-norms.
In doing so, we show that the quasi-norm can be characterized by
using Peetre type maximal function $\vphi_j^*f$.

Let $\Phi, \vphi, \Psi, \psi \in C_0^{\infty}(\R )$ satisfy \cite{FJW,E95}
\begin{align*}
&i)\quad\; \supp\; \Phi, \supp\; \Psi \subset\{ |\xi|\leq 1\}; \;
|\Phi (\xi)|, |\Psi (\xi)| \geq c>0 \; \text{if}\;  |\xi|\leq \frac{1}{2};
\\
&ii)\quad\;
\supp \; \vphi, \supp\; \psi \subset \{\frac{1}{4}\leq |\xi|\leq
1 \}; \;
 |\vphi ( \xi ) |, |\psi ( \xi )|\geq c>0\;  \mathrm{if}\;
 \frac{3}{8} \leq |\xi|\leq
 \frac{7}{8};\\
&iii) \quad\; \Phi (\xi) \Psi (\xi) + \sum_{j=1}^{\infty} \vphi(2^{-j}\xi)
\psi (2^{-j} \xi)=1,  \;\;\forall \;\xi \in \R\, .
\end{align*}

\begin{definition}
The Triebel-Lizorkin space
${F}_p^{\alpha,q}:= {F}_p^{ \alpha,q } (H)
$ associated with $H$ is the completion of the set $\{f\in L^2:
\Vert f\Vert_{F_p^{ \alpha,q } } <\infty\}$ with respect to the
norm $\Vert \cdot \Vert_{F_p^{\alpha,q}}$,
where
$$\Vert f\Vert_{ {F}_p^{\alpha,q} } = \Vert (\sum_{j=0}^{\infty}
2^{j\alpha q}|\vphi_j(H)f|^q )^{1/q}\Vert_p. $$
\end{definition}

Note, that  $1<p<\infty, q=2$, $F^{\al,q}_p(H)$ corresponds to
Sobolev type spaces \cite[p.15]{Tr92}.

For $j\in\Z$
 let $\phi_j(x)=\phi(2^{-j}x)$ and
$K_j(x,y)=\phi( 2^{-j} H ) (x,y)$.

As in \cite{BZ04},
the following lemma is essential to establish the Peetre type maximal
inequalities Lemma \ref{lem:der-phi*} and Lemma \ref{lem:phi*-M}.
We postpone the proof till Section \ref{S4}. To simplify the notation
we set
\begin{equation}\label{eq-wj}
w_j(x):=(1+2^{j/2}|x|)\, .
\end{equation}
\begin{lemma}\label{lem:ker-decay}
For each $n\in \mathbb{N}_0$ there exist
constants $C_n, D_n>0$ such that
\begin{enumerate}
\item
$\displaystyle{
\qquad \vert K_j(x,y) \vert \leq C_n  2^{j/2}w_j(x-y)^{-n} ;}$
\item\label{Th3.2a}
$\displaystyle{ \qquad \vert \frac{\partial}{\partial x}K_j(x,y)\vert
\leq D_n2^{j}w_j(x-y)^{-n}\, .}$
\end{enumerate}
\end{lemma}

Note that the statement here is \textit{simpler}
than the corresponding statements in Lemmas 3.1-3.3
 in \cite{BZ04}. 
\medskip

\begin{remark}
If $ H_\alpha = -\triangle +V_\alpha $ on $\R^d$, where $ V_\alpha (x) = \alpha^2
V(\alpha x )$ is the regular scaling, then for $\phi \in L^\infty$
we have:
\begin{equation}\label{eq:ker-scale}
 \phi (H_\alpha)(x,y) = \alpha^d  \phi(H)(\alpha x,\alpha y ).
\end{equation}

To see this, we note that $\delta_\al \circ \,H\circ\,
\delta_{\al^{-1} } = \al^{-2} H_\alpha  $, hence $\delta_\al \phi (H)
\delta_{\al^{-1} } = \phi( \al^{-2} H_\alpha ) $.
>From this we obtain
$$
\phi (H)(x,y) = \alpha^{-d } \phi(H_\al) (\alpha^{-1} x,\alpha^{-1} y ).
$$

If $ H_h = -\triangle +V_h, V_h= V(\cdot -h) $ is a translation of $V,
h\in \R^d$, we have for $\phi \in L^\infty$,
\begin{equation}\label{eq:ker-trans}
 \phi (H_h )(x,y) =   \phi(H)( x-h, y-h ),
\end{equation}
where we note that $H_h= \tau_h \circ H \circ \tau_{-h} $, hence $\phi(H_h)
= \tau_h \circ \phi(H)\circ \tau_{-h} $.

The above identities mean that the decay estimates for the kernels, hence the
results in this paper, still hold true  for the case when our $V$  is
replaced by its
regular scaling  or translation.
\end{remark}
\medskip

For $s>0$ define the \emph{Peetre maximal function} for $H$ by:
\begin{equation}\label{eq:phi*f}
\phi_j^*f(x):=\phi_{j,s}^*f (x) = \sup_{t\in \R} \frac{|\phi_j(H)f(t)|}{w_j(x-t)^s}
\end{equation}
\noindent
and
$$
\phi_j^{**} f (x):=\phi_{j,s}^{**} f (x) = \sup_{t\in \R}
\frac{| (\phi_j(H)f)'(t)|}{w_j(x-t)^s}.$$

\begin{lemma} \label{lem:der-phi*}
Let $\phi\in C_0^\infty$ and $s>0$. Then
there exists a constant $c_s>0$ such  that
$$
\phi_j^{**}f(x) \leq C_s 2^{j/2} \\*
  \phi_j^{*} f( x).$$
\end{lemma}

\begin{proof}
Using the diadic partition of unity from the beginning
of the section and noting that the sum contains at most
three non-zero terms we write using the convention
that $\psi_{-1}=0$:
\begin{eqnarray*}
\frac{d}{dt}( \phi_j(H)f)(t)& =& \frac{d}{dt} \sum_{\nu=-1}^1
(\vphi \psi)_{j+\nu}(H) \phi_j(H)f(t)\\
&=&  \sum_{\nu = -1}^1 \int_\R \frac{\partial}{\partial t}K_{j+\nu}(t,y)
\phi_j(H)f(y)\, dy\, ,
\end{eqnarray*}
where $K_{j+\nu}(t,y)$ denotes the kernel of $(\vphi \psi)_{j+\nu}(H)$.

Apply Lemma \ref{lem:ker-decay} to obtain
$$
\frac{ |\frac{d}{dt}( \phi_j(H)f)(t)|}{w_j(x-t)^s} \leq
C_n
 \sum_{\nu=-1}^1
2^{j+\nu}
\frac{ |\phi_j(H)f(y)|}{w_{j+\nu}(t-y)^n w_j(x-t)^s}\, dy
\, .$$
We now claim that
\begin{equation}\label{eq:claim}
\frac{ |\phi_j(H)f)(y)| }{w_j(x-t)^s}
\leq \max\limits_{k} \phi_j^*f(x)w_j(t-y)^s\, .
\end{equation}
To confirm the claim we note that
the left hand side of (\ref{eq:claim}) is bounded by
\begin{eqnarray*} \frac{ |\phi_j(H)f)(y)| }{w_j(x-y)^s} \cdot
\frac{w_j(x-y)^s}
{w_j(x-t)^s}
&\le&
\phi_j^*f(x )\frac{w_j(x-t)^s
w_j(y-t)^s} {w_j(x-t)^s}\\
&=&  \phi_j^*f(x)w_j(t-y)^s\, ,
\end{eqnarray*}
It follows that
\begin{eqnarray*}
\frac{ |\frac{d}{dt}( \phi_j(H)f)(t)|}{w_j(x-t)^s}
&\leq &C_n \sum_{\nu=-1}^1 2^{j+\nu}
\phi_j^*f(x )
\int_\R \frac{w_j(t-y)^s}{w_{j+\nu}(t-y)^s}\,  dy \\
&\leq& C_{n} 2^{(j+n)/2}  \phi_j^*f(x)
 \int_\R
\frac{1}{w_j(t-y)^{n-s}}\,
dy\\
&\leq& C_{n,s} 2^{j/2}
\phi_j^*f(x )\,
\end{eqnarray*}
provided $n-s>1 \leftrightarrow n=[s+1]+1$.
This proves Lemma \ref{lem:der-phi*}.
\end{proof}

We are ready to show Peetre maximal inequality.

\begin{lemma}\label{lem:phi*-M} Denote by $M$ the
Hardy-Littlewood maximal operator.
If $0<r<\infty$ there exists a constant $C>0$
such that
$$
 \phi_j^*f(x)
C_{1/r} [M({|\phi_j(H)f|}^r)]^{1/r}(x)\, .$$
\end{lemma}
\begin{proof} The same as those found in
\cite{BZ04,E95}, which is a  modification of the proof in the Fourier case
 \cite[p.16]{Tr83}.
\end{proof}
\smallskip

\begin{remark}
It is known that the Hardy-Littlewood maximal operator
$M$ is bounded on $L^p, 1<p<\infty$.  Lemma \ref{lem:phi*-M} implies that
\begin{equation}\label{eq:phi*<phi}
\Vert  \phi_j^*f \Vert_p \leq c \Vert  \phi_j(H)f \Vert_p, \;\; 0<p \le \infty
\end{equation}
Moreover, if $1<p,q<\infty$ and $\{f_j\}$ is a sequence of functions, then
\begin{equation}\label{eq:F-S-max}
 \Vert  (\sum_j |M f_j|^q)^{1/q}\Vert_{L^p}
\leq C_{p,q}
\Vert (\sum_j |f_j|^q)^{1/q}  \Vert_{L^p},
\end{equation}
by the Fefferman-Stein vector-valued
H-L maximal function inequality.
\end{remark}
\smallskip

We are ready to state the main theorem on
characterization of $ F_p^{\alpha,q}(H)$
using the Peetre maximal function.

\begin{theorem} \label{th:phi*F-inhomo} Let $\alpha\in \R,  0<p,q \leq \infty$.
If $\vphi_j^*f$
is defined for $j\geq 0 $ with $s>1/\min(p,q)$, we have for $f\in L^2 (\R)$
\begin{equation}
 \Vert f\Vert_{F_p^{\alpha,q}} \sim \Vert {\Phi}^*f \Vert_p+
\Vert
\left(\sum_{j=1}^{\infty} ( 2^{j\alpha} \vphi_j^*f)^q \right)^{1/q} \Vert_p.
\end{equation}
Furthermore, $F_p^{\alpha,q}$ is a quasi-Banach space (Banach space if
$p\geq 1, q \geq 1$) and it is independent of the choice of $\{\Phi,
{\vphi}_j\}_{j\geq 1}$.
\end{theorem}

\begin{proof}
Noting that $\phi_j^*f(x) \geq |\phi_j(H)f(x)|$,
it is sufficient to show that
$$ \Vert \Phi^*f\Vert_p +
\Vert (\sum_{j=1}^{\infty}(2^{j \alpha } \vphi_j^* f)^q)^{1/q}  \Vert_p\leq C
\Vert f \Vert_{ F_p^{ \alpha , q } }
$$
but this follows from (\ref{eq:F-S-max}) in the following way:
\begin{eqnarray*}
\Vert \{2^{j\alpha}  \vphi_j^*f\} \Vert_{L^p(\ell^q)}
&\leq & C_{\epsilon,s}
\Vert \{2^{j\alpha}[ M (|\vphi_j(H)f|^r) ]^{1/r} \}\Vert_{L^p(\ell^q)}\\
&=& C_s \Vert \{ \sum_0^{\infty} [M( 2^{j\alpha r}
|\vphi_j(H)f|^r)]^{q/r} \}^{ r/q} {\Vert}^{1/r}_{L^{p/r} } \\
&\leq& C_{s,p,q}
\Vert \{ 2^{j\alpha} \vphi_j(H)f \} \Vert_{L^p(\ell^q)}
\\
&=&C_{s,p,q}\Vert f\Vert_{F_p^{\alpha,q}}\, .
\end{eqnarray*}

Next we show that $F_p^{\alpha,q}$ is independent of the generating
function $\phi $, i.e. if $\tilde{\phi}$ is another function
given a diadic partition of unity then for $f\in L^2$,
$\Vert f\Vert^{\phi}_{F_p^{\alpha,q}}$ and $\Vert
f\Vert^{\tilde{\phi}}_{F_p^{\alpha , q } }$ are
equivalent quasi-norms on $F_p^{\alpha,q}$.

Write $\phi_j(H)= \sum^1_{\nu =- 1}\phi_j(H)\tilde{\phi}_{j+\nu}(H)$
by the identity $\phi_j(x)=\phi_j(x) \sum^1_{\nu =-1}\tilde{\phi}_{j+\nu}(x)$,
for $x\in \R, j\geq 0$, where as before we use the convention that
$\psi_{-1} \equiv 0 $.

We have by Lemma \ref{lem:ker-decay} (1),
\begin{eqnarray*}
\vert \phi_j(H)f(x) \vert &\leq &2^{j/2}\sum^1_{\nu =- 1}
 \int_\R \frac{\vert\tilde{\phi}_{j+\nu}(H)f(y)\vert}{w_j(x-y)^n}
\, dy\\
&\leq& C\sum^1_{\nu =- 1}  2^{j/2} \tilde{\phi}_{j+\nu}^*f(x) \int_\R
\frac{w_j(x-y)^s}{w_j(x-y)^n}\, dy\\
&=& C_N \sum^1_{\nu =- 1}  \tilde{\phi}_{j+\nu}^* f(x),
\end{eqnarray*}
provided $n-s>1$. Thus, for $f\in L^2$
\begin{equation*}
\Vert \{2^{j\alpha}  \vphi_j(H)f\} \Vert_{L^p(\ell^q)}
\leq C_{s,p,q}
\Vert \{ 2^{j\alpha} \tilde{\phi}_j^*f \} \Vert_{L^p(\ell^q)}
\sim \Vert f\Vert^\psi_{F_p^{\alpha,q}},
\end{equation*}
This concludes the proof of Theorem \ref{th:phi*F-inhomo}. 
\end{proof}

As expected from Lemma \ref{lem:phi*-M} we can define the homogeneous
Triebel-Lizorkin spaces
$\dot{F}_p^{\alpha , q } $ and obtain a maximal function
characterization as well.

\begin{definition}\label{d-Fpqa-homo} The  homogeneous Triebel-Lizorkin space
${\dot{F}}_p^{\alpha,q}:= {\dot{F}}_p^{ \alpha,q } (H)
$ associated with $H$ is the completion of the set $\{f\in L^2:
\Vert f\Vert_{ {\dot{F}}_p^{ \alpha,q } } <\infty\}$ with respect to the
norm $\Vert \cdot \Vert_{{ \dot{F}}_p^{\alpha,q}}$,
where
$$\Vert f\Vert_{ {\dot{F}}_p^{\alpha,q} } = \Vert (\sum_{j=-\infty}^{\infty}
(2^{j\alpha}\phi_j(H)f )^q )^{1/q}\Vert_p. $$
\end{definition}

\begin{theorem}\label{th:phi*F-homo}
Let $\alpha\in \R,  0<p,q \leq \infty$. If $\phi_j^*f$
is defined for $j\in \Z $ with $ s>1/\min(p,q) $, then for $f\in L^2$
$$
 \Vert f\Vert_{ {\dot {F} }_p^{ \alpha , q }} \sim \Vert
\left(\sum_{j=-\infty}^{\infty} ( 2^{j\alpha} \vphi_j^*f
)^q \right)^{1/q}\Vert_p\, .$$
Furthermore, $\Vert \cdot \Vert^{\vphi}_{\dot{F}_p^{\alpha,q}}$ and
$\Vert \cdot \Vert^{\psi}_{\dot{F}_p^{\alpha,q}}$ are  equivalent norms
on the  quasi-Banach space $\dot{B}_p^{\alpha,q } $.
\end{theorem}

The proof is completely implicit in that of Theorem \ref{th:phi*F-inhomo} and hence
omitted.

Moreover, like in the Fourier case and \cite{Tr83}, Peetre
maximal inequality allows us to define and characterize
$H$-Besov spaces.

Let $\alpha\in \R,  0<p< \infty, q \leq \infty$.
We define the \emph{Besov space} $B_p^{\alpha,q}:=
{B}_p^{\alpha,q}(H)$ associated with $H$ to be the completion of the
set $\{f\in L^2:
\Vert f\Vert_{B_p^{\alpha,q} } <\infty\}$ with respect to the
norm $\Vert \cdot \Vert_{B_p^{\alpha,q}}$, where

$$ \Vert f\Vert_{B_p^{\alpha,q} } = \Vert \Phi (H)f \Vert_p +
\{\sum_{j=1}^{\infty}
 2^{j\alpha q} \Vert \vphi_j(H)f \Vert_{L^p}^q \}^{1/q}.
$$
Thus $B_p^{\alpha,q} $ is a quasi Banach space (Banach space if $p,q
\geq 1$). The proof of the following Theorem is analogous to
that of Theorem \ref{th:phi*F-inhomo} and thus left to the reader.

\begin{theorem}\label{th:phi*B-inhomo}  Let $\alpha\in \R,
 0<p<\infty, 0<q \leq \infty$. If
$\Phi^*f, \;\phi_j^*f$
is defined for $j>0 $ with $s>1/p$, then for $f\in L^2$
$$
 \Vert f\Vert_{ B_p^{\alpha,q}} \sim
\Vert \Phi^*f\Vert_p + (\sum_{j=1}^{\infty}  2^{j\alpha q}
\Vert \vphi_j^*f\Vert_{L^p}^q)^{1/q}.$$
Furthermore, $B_p^{\alpha,q}$ is well defined and independent of the
choice of generating functions $\vphi=\vphi_0$.
\end{theorem}

Choosing $\eps$ small enough in Lemma \ref{lem:phi*-M} we obtain by
taking $r<p(\leftrightarrow s>1/p)$
\begin{eqnarray*}
\Vert \{2^{j\alpha}\phi_j^*f\}\Vert_{B_p^{\alpha,q}}
& \leq & C_{s} (\sum_j 2^{j\alpha q}
\Vert [M(|\phi_j(H)f|^r)]^{1/r}\Vert_{L^{p}}^q)^{1/q}\\
&=& C_{s} (\sum_j 2^{j\alpha q}
\Vert M(|\phi_j(H)f|^r)\Vert_{L^{p/r}}^{q/r})^{1/q}\\
&\le&C(\sum_j 2^{j\alpha q}\Vert |\phi_j(H)f|^r\Vert_{L^{p/r}}^{q/r})^{1/q}
\\
&=&C (\sum_j 2^{j\alpha q}\Vert\phi_j(H)f\Vert_p^q)^{1/q}
\\
&=& C_{s} \Vert f\Vert_{B_p^{\alpha,q}},
\end{eqnarray*}
where we used the
Hardy-Littlewood maximal function inequality.

The definition and characterization for the homogeneous
Besov spaces $\dot{B}_p^{\alpha,q }$ can be formulated similarly.

\section{High energy estimates}\label{S4}
In this section we proof Lemma \ref{lem:ker-decay}.
We divide the  estimates into high and low energy cases.

Recall from Theorem \ref{th-exk} that
$$
e(x,k)=(\mathrm{sign}(k))^{n}\left(\prod_{j=1}^n\frac{1}{j+i | k|}
\right)\,P_n(x,ik)e^{ikx}$$
where $P_n(x,k)=p_n(\th (x),ik)$ is defined by the recursion formula
$$
p_n (\th x, ik)=
\frac{d\, }{dx}p_{n-1}(\th x, ik)+( ik-n\th x) p_{n-1}(\th x, ik)\, .$$
By the same Theorem
the function
$$e(x,k)\overline{e(y,k)}=
\left(\prod_{j=1}^n\frac{1}{j^2+ k^2}
\right)\, P_n(x,ik)P_n(-y,ik) e^{ik(x-y)}$$
is analytic in $k\in
\R$ and $e (x,-k ) = e (-x,k )$. Furthermore we will
also use the obviously
relation
$\phi(\lambda^2 H)(x,y)  = \overline{\bar{\phi}(\lambda^2 H)(y,x) }$.
For simplicity we write
\begin{equation}\label{def-Qxyk}
Q(x,y,k)=P(x,ik)P(y,-ik)=P(x,ik)\overline{P(y,ik)}\, .
\end{equation}

In this section we set $\lambda=2^{-j/2}$
and $\psi (x) = \phi (x^2)$, where
$ \supp \phi\subset [\frac{1}{4}, 1]$ is fixed. Then
\begin{eqnarray*}
K_j(x,y)= \frac{1}{2\pi} \int \psi (\lambda
\xi ) e(x, k)\overline{e(y,k )}\, dk.
\end{eqnarray*}

\subsection{High and low energy estimates} Using the above explicit formula
for $e(x,k)\overline{e(x,k)}$ we set
$$I^+(x,y)=\int_0^\infty
\phi_j(k^2)\left(\prod_{j=1}^n \frac{1}{j^2+k^2}\right)
Q(x,y,k) \, e^{ik(x-y)}dk$$
and
$$I^-(x,y)=\int_{-\infty}^0
\phi_j(k^2) \left(\prod_{j=1}^n\frac{1}{j^2+k^2}\right)
Q(x,y,k)\,  e^{ik(x-y)}dk\, .$$
Then we have
\begin{equation}\label{eq-phiH}
2\pi \phi_j(H)\vert_{H_c} (x,y) =I^+(x,y)+I^-(x,y)\, .
\end{equation}

It is enough to deal  $I^+(x,y)$ since $I^-(x,y)=I^+(-x,-y)$
by a change of variable $k\rightarrow -k$.

We note, that
$$k\mapsto R(x,y,k)=R(k)= \left(\prod_{j=1}^n \frac{1}{j^2+k^2}\right)
 Q(x,y,k)$$
is a rational function with numerator and denumerator both
of degree $2n$, and uniformly bounded it $x$ and $y$.
It follows that for $|k|>\epsilon >0$, we have
$|\frac{d^j\, }{dk^j}R(x,y,k)|\le C_j |k|^{-j}$
where $C_j>0$. Furthermore
$|\frac{d^j\, }{dk^j}\psi (\lambda k)|
\le \|\psi^{(n)}\|_\infty \lambda^j$.
Note, that the left hand in only non-zero for
$|k|\in [2^{j/2-1},2^{j/2}]$. This
gives first of all an uniform estimate, independent of
$x$ and $y$:
$$|I^+(x,y)|\le \Vert R\Vert_\infty \int_{2^{j/2-1}}^{2^{j/2}}|\psi (\lambda k)|\, dk
\le 2^{-1}\|R\|_\infty\|\psi\|_\infty \lambda^{-1}\, .$$
Let $I_j=[2^{j/2-1},2^{j/2}]$ and note that
$|I_j|=2^{j/2-1}$, where $|\cdot |$ stand for the
standard Lebesgue measure. Then integration by parts gives:
\begin{eqnarray*}
\vert I^+ (x,y)\vert &=&\left|
\frac{(-1)^n}{i^n(x-y)^n} \int_0^\infty  \frac{d^n}{d k^n}[\psi (\lam
 k) R(k)]e^{ik(x-y)}\, dk \right| \\
&\leq &|x-y|^{-n} \int_0^\infty \left| \sum_{s=0}^n
\left(\begin{array}{c}n\\
s\end{array}\right)\frac{d^s\, }{d k^s }(\psi(\lam
k)) \frac{d^{n-s}\, }{dk^{n-s}}R(k)  \right|\,  dk\\
&\le & C_n |x-y|^{-n} \sum_{s=0}^n\lambda^s\left(\begin{array}{c}n\\
s\end{array}\right)
\int_{|k|\in I_j}
\left| \frac{d^{n-s}\, }{dk^{n-s}}R(k)\right|\,  dk\\
&=& C'_n \lam^{n-1}/ {|x-y|^{n} },
\end{eqnarray*}
where  we used
that
$$\lambda^s \int_{|k|\in I_j}
\left| \frac{d^{n-s}\, }{dk^{n-s}}R(k)\right|\,  dk
\le C\lambda^{n-1}$$
for all $s=0,1,\ldots ,n$.
It follows, that there exists a constant $C_n$ such that
\begin{equation}\vert I^+ (x,y)\vert \le  C_n 2^{j/2}/ (1+2^{j/2} |x-y|)^{n}
=C_n2^{j/2}w_j(x-y)^{-n}\, .
\end{equation}

The same estimate holds also for for  $I^- (x,y)$.
Hence it it follows that
\begin{equation} \vert \phi_j(H) |_{ {\cal H}_c} (x,y)\vert \le
C_n \lam^{-1}/ (1+\lam^{-1} |x-y|)^{n}\,   \qquad j>0\, .
\end{equation}

The low energy estimates follows in the same way
by setting $j=0$ using that
$$
2\pi \Phi(H) |_{ {\cal H}_c} (x,y)
= \int_{-1}^{1}\Phi(k^2) |R (x,y,k)| e^{ik(x-y)}dk.
$$
Hence
$$ \vert \Phi(H) |_{ {\cal H}_c} (x,y)\vert \le
C_n  (1+|x-y|)^{-n}\, .$$

\subsection{Proof of Lemma \ref{lem:ker-decay} part 2}
In the following the subscript
${\cal H_c}$ is suppressed, and we denote
by $\phi_j(H) (x,y) $
the kernel of the continuous part of the operator $\phi_j(H) (x,y) $.

According to equation \ref{eq:kernel-repre}, we know
the kernel $\phi_j(H)(x,y)=K_j(x,y)$ ($j\ge 0$)  has a unified expression:
$$\phi_j(H)(x,y)=\frac{1}{2\pi}\int \phi_j(k^2)R(x,y,k)e^{ik(x-y)}\, dk\, .$$
As $\phi_j(k^2)$ has compact support, and
$x\mapsto \frac{d}{dx}\phi_j(k^2)R(x,y,k)e^{ik(x-y)}$ is smooth and hence bounded
on the support of $\phi_j(k^2)$, we can interchange the order
of integration and differentiation to get:
\begin{eqnarray*}
 2\pi \frac{\partial}{\partial x}  \phi_j(H) (x,y)& =&
\frac{\partial\, }{\partial x}\int \phi_j(k^2)
R(x,y,k)e^{ik(x-y)}\, dk\\
&=&\int \phi_j(k^2)\frac{\partial\,}{\partial x}[
R(x,y,k) e^{ik(x-y)}]\, dk,\\
&=& \int \phi_j(k^2)|A(k)|^2[ikP(x,ik)+\frac{\partial\, }{\partial x}P(x,ik)]
P(y,-ik)e^{ik(x-y)}\, dk
\end{eqnarray*}
The function $[\partial/\partial x]P(x,ik)$
is a polynomial of $\th x$ and $ik$ having degree $n$ in $k$.
Note also, that
$$  \left|\frac{d^i}{dk^i} [k \phi_j(k^2)]\right| =\lambda^i
C(k)\Vert \phi_j^{(i)}\Vert_\infty=O(\lambda^{i-1} )$$
where $C(k)$ is a polynomial in $k$, $i\ge 0$.
By the same arguments as in the beginning of this
section, we get
$$
\left|  \frac{\partial}{\partial x} \phi_j(H) (x,y) \right| \le C_n \lam^{n-2}
 (1+|x-y|)^{-n},  \qquad n\ge 1.$$
Together with the fact that
$$
\left|  \frac{\partial}{\partial x} \phi_j(H) (x,y) \right| \le O (\lam^{-2})
$$
this implies that
$$
\left|  \frac{\partial}{\partial y} \phi_j(H) (x,y) \right| \le C_n \lam^{-2}
 (1+\lam^{-1}|x-y|)^{-n}$$
completing the proof of Lemma \ref{lem:ker-decay} (2).
\hfill $\Box$

\section{Interpolation of $F_p^{0,2} \sim L^p,\; 1<p<\infty$ }
\noindent
Let $\phi_j,\psi_j$ be as in $\S$3. We may assume that
$\Vert\Phi\Vert_\infty$, $\Vert\phi\Vert_\infty$,
$\Vert\Psi\Vert_\infty$, $\Vert\psi\Vert_\infty$ are all $\le 1 $.
Let $\phi_0=\Phi$. Define
$$ Q: L^2 \rightarrow L^2(\ell^2) \quad
Qf = \{\phi_j(H)f\}_0^\infty \, .$$
and
$$ R: L^2(\ell^2)\rightarrow L^2 \quad
R( \{ g_j\}_0^\infty ) = \sum_{j=0}^\infty R_j g_j
=\sum_{j=0}^\infty \psi_j(H)g_j\, .$$

Note that $ \Vert f \Vert_{F_p^{0, 2}(H)} = \Vert Qf \Vert_{L^p(\ell^2)}$

It is easy to see that $ RQ=Id: L^2\rightarrow L^2$ and
$ QR \le 3Id:  L^2(\ell^2)\rightarrow L^2(\ell^2)$.
We will use $Q$ and $R$ to identify $F_p^{0, 2}(H)$ with $L^p$.

\begin{theorem} \label{th:Lp-interp}
Let $1<p<\infty$. Then $F_p^{0,2}$ and $ L^p$ are isomorphic
and have equivalent norms.
\end{theorem}

To prove the theorem, we first show that $Q$: $L^p\to L^p(\ell^2)$ and
$R$: $L^p(\ell^2) \to L^p$, $1<p<\infty$.
First we note that
\begin{equation}\label{eq:F<Lp}
\Vert f\Vert_{F_p^{0,2}(H)}\lesssim \Vert f\Vert_p
\end{equation}
and
\begin{equation}\label{eq:Lp<F}
\Vert g\Vert_p \lesssim \Vert g\Vert_{F_p^{0,2}(H)}
\end{equation}
for $f,g\in L^2\cap L^p\cap F_p^{0,2}(H)$.
Here
(\ref{eq:F<Lp}) follows directly from the definition of $Q$
(and $F_p^{0,2}(H)$) and equation (\ref{eq:Lp<F}) follows from the
identity $R\circ Q=id_{L^2}$, i.e.,
$\sum \vphi(H)\psi_j(H)=id$.
By continuity extension of $Q,R$,  
(\ref{eq:F<Lp}) and (\ref{eq:Lp<F}) hold for all of 
$F_p^{0,2}(H)$ and $L^p$.
 So we obtain the continuous embedding
\begin{align}
&L^p\subset F_p^{0,2}(H)\label{eq:L-F}\\
& F_p^{0,2}(H) \subset L^p,\label{eq:F-L}
\end{align}
which proves Theorem \ref{th:Lp-interp}.

The remaining part of this section will be devoted to
showing the boundedness of $Q$ and $R$. In the following,
Lemma \ref{lem:Q-R-L2}
and Lemma \ref{lem:Q-w-L1} imply that $Q$ is bounded
from $L^p$ to $L^p(\ell^2)$, and,
Lemma \ref{lem:Q-R-L2}
and Lemma \ref{lem:R-w-L1} imply that $R$ is bounded
from $L^p(\ell^2)$ to $L^p$ by interpolation and duality.

\begin{lemma} \label{lem:Q-R-L2} $Q: L^2\rightarrow   L^2(\ell^2)$ and
$R: L^2(\ell^2)\rightarrow L^2$
are well-defined bounded operators.
\end{lemma}
\begin{proof} Let $\{ g_j\}\in  L^2(\ell^2)$.
Note that $R_j= \psi_j(H)$ is bounded on $L^2$:
$\Vert R_j g\Vert_2\le \Vert\psi\Vert_\infty\Vert g\Vert_2$.
Thus
\begin{eqnarray*}
(\sum_{j=0}^\infty R_j g_j, \sum_{j=0}^\infty R_j g_j)
&=&\sum_{\nu=-1}^{1}\sum_{j=0}^\infty (R_j g_j,R_{j+\nu}g_j)\\
&\le&\sum_{\nu=-1}^1\sum_j \Vert R_j g_j\Vert_2\Vert R_{j+\nu}g_j\Vert_2\\
&\le& 3\sum_j\Vert  g_j\Vert_2^2\\
&= &3 \Vert g_j\Vert_{L^2(\ell^2)}^2\, .
\end{eqnarray*}
\end{proof}

We now derive some needed estimates for the kernel of $Q_j=\phi_j(H)$.
As before, $Q_j(x,y)$ denotes the kernel of $\phi_j(H)E_{ac}=
\phi_j(H_{ac})$.

Define $$
 \tilde{Q}_j(x,y)=
\left\{
\begin{array}{ll}
\int_{I_k} Q_j(x,y) b_k (y)dy &if \; 2^{j/2}|I_k|\ge 1\\
\int_{I_k} (Q_j(x,y)- Q_j(x,\bar{y})) b_k (y)dy &if \; 2^{j/2}|I_k|< 1.
\end{array}\right.
$$
\begin{lemma}\label{lem:ker-Qj} Let $ y\in I =(\bar{y}-\frac{t}{2}, \bar{y}+\frac{t}{2})
$, $t= |I| $ and $I^* =(\bar{y}-t, \bar{y}+t)$. Then the
following holds:
\begin{enumerate}
\item  If $2^{j/2 } |I| \ge 1$, then
$$ \sup_{y\in I}\int_{R \setminus I^*} |Q_j(x,y) | dx \le C (2^{j/2 } |I|)^{-1}
$$
\item
If $ 2^{j/2 } |I| < 1$, then
$$\sup_{y\in I}\int_{R\setminus I^*} |Q_j(x, y)-Q_j(x, \bar{y} )  | dx
\le C  2^{j/2 } |I|.
$$
\end{enumerate}
In particular we have
\begin{equation}\label{eq:tilde-Q}
\int_{R\setminus I^*}\sum_j|\tilde{Q}_j(x,y)|\le 3C.
\end{equation}
\end{lemma}

\begin{proof}
For the first inequality, let $2^{j/2 } |I|\ge 1$.
\begin{eqnarray*}
\int_{\R \setminus I^*} |Q_j(x,y) |\,  dx
&\le& C_n\int_{|x-{y}|>t/2}
\frac{2^{j/2}}{(1+2^{j/2}|x-{y}| )^n }\, dx \\
&\le&C (2^{j/2 } |I|)^{-1}, \qquad (n=2).
\end{eqnarray*}

If  $ 2^{j/2 } |I| < 1$,
$y\in I$,  $\bar{y} =$ center of $I =(\bar{y}-\frac{t}{2}, \bar{y}+\frac{t}{2})$.
\begin{eqnarray*}
\int_{\R \setminus I^*} |Q_j(x,y)-Q_j(x,\bar{y} ) |\,  dx
&=& \int_{\R \setminus I^*} |\int_{\bar{y}}^y \frac{\partial}{\partial z
}\psi_j(H) (x,z) dz | dx\\
&\le & C_n |y-\bar{y} |\int_{|x-\bar{y}|>t}
\frac{2^{j}}{(1+2^{j/2-1}|x-\bar{y}| )^n }\,  dx \\
&\le& C_n 2^{j/2} t \int_{|u|>(2^{-j/2+1} )^{-1} t } \frac{ du }{(1+
|u|)^n} \qquad (n=2)\\
&\le&C 2^{j/2 } |I|.
\end{eqnarray*}

Now (\ref{eq:tilde-Q}) follows easily from the above two inequalities.
\end{proof}

\begin{lemma}\label{lem:Q-w-L1} $Q$ is bounded from $L^1$ to weak  $ L^1(\ell^2)$, i.e.,
$$ \vert \{ x: (\sum_0^\infty |Q_j f(x) |^2)^{1/2} >\lam\}\vert
\le C\lam^{-1}\Vert
f\Vert_1,  \quad \lam >0.
$$
\end{lemma}

\begin{proof} Let $f\in L^1\cap L^2$. 
By the Calder\'on-Zygmund decomposition, 
there exists a sequence of intervals $ \{I_k\}$ and
and functions  $\{b_k\}$ with $\supp \;b_k\subset I_k$
such that with $g=\sum_k b_k$ we have
$f=g+b$, for some $g\in L^2$.
Furthermore:

\begin{enumerate}
\item[i)] $|g(x)|\le C \lam $ a.e.

\item[ii)] $b_k(x)=f(x)- |I_k|^{-1} \int_{I_k} f dx
$

\item[iii)] $\lam \le |I_k|^{-1} \int_{I_k} |f| dx \le 2\lam$.

\item[iv)] $ |D_\lam|= |\cup_k I_k|=\sum_k |I_k|
\le \lam^{-1}\Vert f \Vert_1.
$
\end{enumerate}
Let $\{e_n\}_m$ be an orthonormal basis in ${\cal H}_{pp}$
and $\lam_m$ the corresponding eigenvalues in $\sigma_p$. Note
that this set is finite.
\begin{eqnarray*}
\int\sum_j |Q_j g(x)|^2 dx &=&\sum_j \Vert \phi_j(H) g \Vert_2^2\\
&=& \sum_j (\Vert\cal{F}^*\phi_j(k^2)\cal{F}g\Vert_2^2 +
\sum_{\sigma_p} \Vert\phi_j(\lam_n)(g,e_m)e_m\Vert_2^2)\\
&=& \sum_j (\Vert\phi_j(k^2)\cal{F}g(k)\Vert_2^2 +
\sum_{\sigma_p} |\phi_j(\lam_n)|^2|(g,e_m)|^2 ) \\
&\le&\int(\sum_j \vert\phi_j(k^2)\vert^2)\vert\cal{F}g(k)\vert^2\,  dk
+ 3\Vert\phi\Vert_\infty^2 \sum_{\sigma_p} \vert (g,e_m)\vert^2 \\
&\le& 3\Vert\phi\Vert_\infty^2\int\vert\cal{F}g(k)\vert^2 dk
+ 3\Vert\phi\Vert_\infty^2 \Vert E_{pp}g\Vert^2 \\
&\le& 3\Vert\phi\Vert_\infty^2\Vert E_{ac}g\Vert^2
+ 3\Vert\phi\Vert_\infty^2 \Vert E_{pp}g\Vert^2 \\
&=& 3\Vert\phi\Vert_\infty^2 \Vert g\Vert^2,
\end{eqnarray*}
where we note that $\cal{F}$: $L^2\to L^2$ surjection and
 $\cal{F}$: $\cal{H}_{ac}\to L^2$ isometry. By
Chebychev inequality we have
$$ \vert \{ x: (\sum_0^\infty |Q_j g(x) |^2)^{1/2} >\lam/2\}\vert
\le C\lam^{-2}\Vert g\Vert_2^2\le C\lam^{-1} \Vert f\Vert_1.
$$
We now only need to show
$$
\vert \{x \notin \cup I_k^*:
(\sum_j  |Q_j b (x)|^2)^{1/2} >\lam/2\}\vert  \le C \lam^{-1}
\Vert f\Vert_1.$$

The left hand side  is bounded by
$$
\frac{2}{\lam } \sum_k \int_{\R\setminus \cup I_k^*}
(\sum_j  |Q_j b_k (x)|^2)^{1/2} dx.
$$
For each $k$, since $\int b_k=0$,
\begin{eqnarray*}
(\sum_j |Q_j b_k (x)|^2)^{1/2}&\le&
\sum_j |Q_j(E_{ac}+E_{pp}) b_k (x)|\\
&\le& \sum_j |\int \tilde{Q}_j(x,y) b_k (y)\, dy|+
\sum_j |\sum_{\sigma_{pp}}\phi_j(\lam_m)(b_k,e_m)e_m(x)|\\
&\le& \int \sum_j |\tilde{Q}_j(x,y) b_k (y)|\, dy+
\sum_{\sigma_{pp}}(\sum_j |\phi_j(\lam_m)|\Vert b_k\Vert_1 |e_m(x)|\\
&\le& \int \sum_j |\tilde{Q}_j(x,y) b_k (y)|\, dy+
3\Vert \phi\Vert_\infty \Vert b_k\Vert_1 \sum_{\sigma_{pp}}|e_m(x)|
\end{eqnarray*}
where we note that $e_m\in L^1\cap L^\infty$
and  $\sigma_{pp}$ is finite, having $n$ eigenvalues.

Equation (\ref{eq:tilde-Q}) gives now the following estimates for the
left hand side
\begin{eqnarray*}
l.h.s. &\le&  \frac{2}{\lam } \sum_k \int_{\R\setminus \cup I_k^*}
\sum_j  \int |\tilde{Q}_j(x,y)| |b_k (y)|\, dy dx+\\
&& 3 \Vert \phi\Vert_\infty \frac{2}{\lam }
\sum_k \int_{\R\setminus \cup I_k^*}\Vert b_k\Vert_1 (\sum_{\sigma_{pp}}|e_m(x)|)\, dx\\
& \le& \frac{2}{\lam } \sum_k\int_{y\in I_k} |b_k (y)|dy\int_{\R\setminus I_k^*}
 \sum_j  |\tilde{Q}_j(x,y)|  dx+  \frac{C}{\lam } \sum_k \Vert b_k \Vert_1
\sum_{\sigma_{pp}}\Vert e_m \Vert_1\\
&\le&  \frac{C'}{\lam } \sum_k\int_{y\in I_k} |b_k (y)|\, dy
\le C\lam^{-1}\Vert f\Vert_1\, .
\end{eqnarray*}
\end{proof}

By Lemma \ref{lem:ker-Qj}, we have the following
estimates which will be used in the proof of Lemma \ref{lem:R-w-L1}:

If $2^{j/2} |I_k |\ge 1$, then
$$ \sum_{j: 2^{j/2} |I_k |\ge 1 } \int_{\R - I^*} \vert R_j(x,y)\vert
\le C |I_k|^{-1} \sum_{2^{j/2} |I_k |\ge 1 } 2^{-j/2} \le C.
$$

If $ 2^{j/2} |I_k |< 1$, then
$$ \sum_{j: 2^{j/2} |I_k| < 1 } \int_{\R - I^*} \vert R_j(x,y)-R_j(x,\bar{y})\vert
\le C |I_k| \sum_{2^{j/2} |I_k |< 1 } 2^{j/2} \le 2C.
$$

\begin{remark} $  R_jf(x)=  
\int  \psi_j(H)(x,y)f(y) dy+
\sum_{\sigma_{pp}} \psi_j(\lam_m)(f,e_m)e_m(x)  $ is
well defined for $f\in L^1\cap L^2$ and for $f\in L^1$ by all appropriate means,
since $\psi_j(H)(x,\cdot) \in L^1\cap L^\infty$.
\end{remark}

\begin{lemma} \label{lem:R-w-L1} Let $R_j= \psi_j(H) $. Then
$R= \{R_j\}$ are bounded from
$ L^1(\ell^2)$ to weak-$L^1$.
\end{lemma}

\begin{proof} We will show that 
there exists a positive constant $C $ such that
\begin{equation}\label{eq:partial-sum-Rf}
\vert \{x: |\sum_0^N  R_j f_j (x) | >\lam \} \vert \le C \lam^{-1} \Vert
 \{ f_j\} \Vert_{L^1(\ell^2) },
\end{equation}
for all $N\in \N$, $\{ f_j \} \in  L^1(\ell^2)$ and $ \lam >0$.
Then, the inequality holds for $N=\infty$ and
$\{ f_j \} \in  L^1(\ell^2)\cap  L^2(\ell^2)$ by passing to the limit,
noting that
$\sum_0^\infty  R_j f_j (x)\in L^2$ if $\{ f_j \} \in  L^2(\ell^2)$.
Finally, the lemma follows by the fact that
$L^1(\ell^2)\cap  L^2(\ell^2)$ is dense in $L^1(\ell^2)$.

Alternatively, suppose $\{ f_j \} \in  L^1(\ell^2)$ only.
Observe that (\ref{eq:partial-sum-Rf})
suggests that the partial sum $R_N f=\sum_0^N  R_j f_j $
is a Cauchy sequence in measure, therefore there exists
a unique function, which we define to be $Rf$,
so that $R_N f$ converges to $Rf$ in measure and
$Rf=\sum_{j=0}^\infty R_jf_j$ satisfies (\ref{eq:partial-sum-Rf}).

So, it is sufficient to show (\ref{eq:partial-sum-Rf}).
Calder\'on-Zygmund decomposition: Let $F(x)=(\sum_{j=0}^\infty |f_j
(x)|^2)^{1/2} \in L^1$. There exists a collection of countable
disjoint open intervals $\{ I_k \} $ such that
\begin{enumerate}
\item[a)]$ |F(x)|\le \lam, \;\;\textrm{a.e.} \; x \in \R\backslash \cup_k I_k$,
\item[b)] $\sum_k |I_k|\le \lam^{-1} \Vert F\Vert_1$,
\item[c)] $\lam\le |I_k|^{-1}\int_{I_k } |F(x)|dx \le 2\lam, \quad \forall k$.
\end{enumerate}

Define
$$ g_j(x)=\left\{
\begin{array}{ll}
|I_k|^{-1}\int_{I_k } f_j dy, & x\in I_k\\
f_j(x)   & otherwise
\end{array}\right.
$$

$$
b_j(x)=\left\{
\begin{array}{ll}
f_j-g_j,   & x\in I_k\\
0 & otherwise
\end{array}\right.
$$
Then, if $x\in \R\setminus \cup_k I_k$, $(\sum_{j=0}^{\infty} |g_j(x)|^2
)^{1/2}= ( \sum_{j=0}^{\infty} |f_j(x)|^2)^{1/2}$.
And, if $x\in  I_k$
\begin{eqnarray*}
(\sum_{j=0}^{\infty} |g_j(x)|^2
)^{1/2} &\le& (\sum_{j=0}^{\infty} |I_k|^{-2}\vert \int_{I_k } f_j(y) \,
dy\vert^2 )^{1/2}\\
&\le &|I_k|^{-1} \int_{I_k }(\sum_{j=0}^{\infty} \vert f_j(y)\vert^2
)^{1/2}\, dy\\
& \le& 2\lam
\end{eqnarray*}
by Minkowski inequality.

We also have
\begin{eqnarray*}
\Vert \{g_j(x)\}\Vert_{L^2(\ell^2)}^2 & = &\sum_{k=0}^{\infty} \int_{I_k}
(\sum_j |g_j(x)|^2)\,  dx +\int_{\R \setminus\cup I_k} (\sum_j |g_j(x)|^2)\, dx\\
&\le & (2\lam)^2 \sum_k |I_k| + \lam \int_{\R \setminus\cup I_k}
(\sum_j |f_j|^2)^{1/2}\,  dx\\
& \le & C\lam \Vert F\Vert_1\, .
\end{eqnarray*}

Finally, by Chebychev inequality and Lemma \ref{lem:Q-R-L2}, we get:
\begin{eqnarray*}
\vert \{x: |\sum_0^N  R_j g_j (x) |  >\lam/2 \} \vert &\le & C \lam^{-2} \Vert
\sum_0^N R_j g_j \Vert^2 \\
& \le & C' \lam^{-2} \Vert \{g_j\} \Vert_{L^1(\ell^2) }^2  \le C \lam^{-1}
\Vert F\Vert_1\, .
\end{eqnarray*}

It remains to show
$$\vert \{x \notin \cup I_k^* : |\sum_0^N  R_j b_j (x) | >\lam/2
\}\vert
\le  C \lam^{-1}
\Vert F\Vert_1\, .$$

The left hand side is not exceeding
$ \frac{2}{\lam } \sum_k \int_{\R\setminus\cup I_k^*} |
\sum_{j=0}^N  R_j b_{j,k} (x)|dx$,
where $  b_{j,k}= b_j \chi_{I_k}$.
For each $k$, define
$$
\tilde{R}_j^k (x,y)= \left\{
\begin{array}{ll}
 R_j (x,y) & \;\textrm{if}\; 2^{j/2} |I_k| \ge 1\\
R_j (x,y)- R_j (x,\bar{y}_k ) & \;\textrm{if}\; 2^{j/2} |I_k| < 1
\end{array}\right.
$$
Note that the domain of $ R_j$ has a natural extension from 
$L^1 \cap  L^2$ to $L^1$. 
Then, since $\int b_{j,k}=0$ and 
by Jensen inequality 
$(\sum_{j=0}^N \vert  \tilde{R}_j^k (x,y)\vert^2
)^{\frac{1}{2}}
\le \sum_{j=0}^N \vert  \tilde{R}_j^k (x,y)\vert $
we obtain 
\begin{eqnarray*}
\int_{\R\setminus I_k^*} |\sum_{j=0}^N  R_j b_{j,k} (x)|\, dx
&=& \int_{\R\setminus I_k^*} |\sum_{j=0}^N
\int_{I_k} \tilde{R}_j^k(x,y) b_{j,k} (y)dy |\, dx\\
&\le& \int_{I_k} ( \sum_{j=0}^N |b_{j,k}|^2 (y) )^{\frac{1}{2}}\,  dy
\int_{\R\setminus I_k^*}
(\sum_{j=0}^N \vert  \tilde{R}_j^k (x,y)\vert^2 )^{\frac{1}{2}}\, dx\\
&\le& C \int_{I_k} ( \sum_{j=0}^N |b_{j,k}(y)|^2 )^{\frac{1}{2}}\, dy\\
&\le& C \int_{I_k} ( \sum_{j=0}^N |f_j|^2 )^{\frac{1}{2} }\, dy +C
 \int_{I_k} ( \sum_{j=0}^N |g_j |^2 )^{\frac{1}{2} }\, dy\\
&\le& 2C  \int_{I_k} ( \sum_{j=0}^\infty |f_j|^2 )^{\frac{1}{2}}\, dy,
\end{eqnarray*}
where we used the pointwise estimate by Minkowski inequality
$$( \sum_0^N |g_j(y) |^2 )^{\frac{1}{2} }
\le ( \sum_0^\infty
(|I_k|^{-1}\int_{I_k} |f_j(y)| )^2)^{\frac{1}{2}}
\le |I_k|^{-1} \int_{I_k} ( \sum_0^N |f_j(y)|^2 )^{\frac{1}{2}}\, dy\, .
$$

Hence
$$\vert \{ x\in \R\setminus \cup I_k^*: \vert \sum_0^N  R_j b_{j} (x) \vert >
\lam /2\} \vert
\le \frac{4C}{\lam} \sum_k \int_{I_k}(\sum_j |f_{j}|^2 )^{\frac{1}{2}}\,  dy
\le \frac{4C}{\lam} \Vert ( \sum_j |f_j|^2  )^{\frac{1}{2}}\Vert_1\, .
$$
as desired.

By now,
we have shown that
\begin{eqnarray*}
 \vert \{ x: \vert \sum_0^N R_j f_{j} (x)\vert >
\lam \} \vert
&\le& \vert \{ x: \vert \sum_0^N R_j g_{j} (x)\vert >
\frac{\lam}{2} \} \vert+ \vert \{ x: \vert \sum_0^N R_j b_{j} (x)\vert >
\frac{\lam}{2} \} \vert\\
& \le &C \lam^{-1} \Vert F\Vert_1,
\end{eqnarray*}
where $C$ is independent of $N$.
Passing to limit, since convergence a.e. implies convergence in
measure, we obtain
$$
 \vert \{ x: \vert \sum_0^\infty \vert R_j f_{j} (x)\vert >
\lam \} \vert
\le C \lam^{-1} \Vert F\Vert_1.
$$
\end{proof}

\section{Remark on time decay of the wave function}\label{S-w}
We conclude the paper with  a time decay result on the
wave function $\psi(t,x)=e^{-itH}f$ based on
on the $H$-Besov space method.
Generalized Besov space method has been considered in e.g. \cite{J92},
\cite{JN94}, \cite{JP85}, \cite{BZ04} 
in the study of PDEs related to
perturbation for Schr\"odinger operators.

By \cite[Theorem 7.1]{BZ04} or \cite[Theorem 5.1]{JN94}
we know that
if $V$ is in the Kato class $\cal{K}$
and if $\cal{D}(H^m)=W_p^{2m}$ for
some $m\in \N $,
$1\le p < \infty$, then for $1\le q \le \infty , 0<\al<m,
B^{\al,q}_p(H) =B^{2\al,q}_p(\R^d) $. 
It is easy to see that
if $V$ is $ C^\infty $ with all derivatives bounded, then the domain
condition on $H$ is verified for all $m\in \N$.
Since $V\sim\sech^2 \,x$ has this property,
$B^{\al,q}_p(H) =B^{2\al,q}_p(\R^d) $ for all $\al>0$.
Combining this with
 \cite[Theorem 5.2]{JN94}, we obtain the time decay
of $\psi(t,x)$ on ordinary Besov spaces. 
Let $\langle t\rangle=(1+t^2)^{1/2}$.

\begin{theorem}\label{th:besov-e{-itH}}
Suppose $ 1\le p<\infty$, $1\le q\leq \infty$ and $\alpha>0$.
 Let $\beta =\vert \frac{1}{2}-\frac{1}{p}\vert$.
Then $e^{-itH}$ maps $B_p^{\alpha+2\beta,q}(\R)$ continuously
to $B_p^{\alpha,q}(\R) $.
Moreover, $e^{-itH}$ maps $B_p^{2\beta,q}(\R)$ continuously
to $L^p$.
In both cases the operator norm is less than or equal to $C\langle t\rangle^\beta$. That is,
$$
\Vert\psi(t,x)\Vert_{B_p^{\alpha,q}(\R)}
\lesssim \langle t\rangle^\beta \Vert f \Vert_{B_p^{\alpha+2\beta,q}(\R)}$$
and
$$
\Vert\psi(t,x)\Vert_{L^p}
\lesssim \langle t\rangle^\beta \Vert f \Vert_{B_p^{2\beta,q}(\R)}.$$
\end{theorem}
Compare Proposition 7.3 in \cite{BZ04}.

Finally we give an interpolation result.
We know that  $B^{\alpha, q}_p(H) = B^{2\alpha,q}_p(\R)$, $ \alpha>0$.
In particular, $F^{\alpha, p}_p(H) = F^{2\alpha,p}_p(\R)$.
On the other hand, by Theorem \ref{th:Lp-interp},
 $F^{0, p}_p(H)=L^p=F^{0, p}_p(\R)$. Thus we obtain
Theorem \ref{th:F-interp} using complex interpolation method. 

\begin{theorem}\label{th:F-interp} If $\al>0$, $1<p<\infty$
and $2p/(p+1)<q<2p$, then
$$F^{\alpha, q}_p(H) = F^{2\alpha,q}_p(\R)$$
and if $\al>0$, $1\le p<\infty$, $1\le q\le \infty$, then
$$B^{\alpha, q}_p(H) = B^{2\alpha,q}_p(\R).$$
\end{theorem}

In general $e^{-itH}$ is not bounded on $L^p$ unless $p=2$.
[JN94,RS04]. However one can obtain the boundedness
on $F$ and $B$ spaces.

Let $S_p$ be a set of numbers $q$ such that
$$
\begin{array}{lcl}
p\le q<2 &\textrm{if}& p<2; \\
2        & \textrm{if}& p=2;\\
2<q\le p &\textrm{if}& p>2.
\end{array}
$$
Based on the mapping property: if $\al>0$,
$1\le p<\infty, 1\le q\le\infty$,
$\beta=|\frac{1}{p}-\frac{1}{2}|$,
$e^{-itH}$: $F^{\al+2\beta,p}_p\to F^{\al,p}_p$
and $F^{2\beta,p}_p\to L^p$,
and the
interpolation $(F^{\al+2\beta,p}_p,F^{2\beta,p}_p)_\theta
\to (F^{\al,p}_p,F^{0,2}_p)_\theta$,
we obtain
\begin{proposition} If $\al_1>0$, $1<p<\infty$,
$q_1\in S_p$, $\beta:=\beta(p)=|\frac{1}{p}-\frac{1}{2}|$,
then
\begin{equation}\label{eq:e-F-bound}
e^{-itH}: \;F^{\al_1+2\beta,p}_p(H)\to F^{\al_1,q_1}_p(H).
\end{equation}
If $\al>0$, $1\le p<\infty, 1\le q_1\le \infty$,
$$B^{\al+2\beta,q}_p(H)\to B^{\al,q}_p(H)$$
with operator bound $C\langle t\rangle^\beta$.
\end{proposition}

If we try to use interpolation coupled with $L^2\to L^2$,
we proceed
 $$(F^{\al_1+2\beta,p}_{p},F^{0,2}_2)_\theta
\to (F^{\al_1,q_1}_p,F^{0,2}_2)_\theta. $$
 We find that for each $t\in (0,1)$
$$
F^{\al_2+2\beta_t,p_t}_{p_t}(H)\to F^{\al_2,q_2}_{p_t}(H)
\quad \;\forall \al_2>0, 1<p_t<\infty$$
where $1/p_t=(1-t)/p+t/2$, $q_2\in S_{p_t}$.
It means the interpolation with $F^{0,2}_2\to F^{0,2}_2$
does not help improve (\ref{eq:e-F-bound}) at all.
This can also be observed from the two-triangle shaped
region bounded by
$\{(\frac{1}{p},\frac{1}{q_1}): q_1\in S_p\}$.

The result for $B$ spaces follows from 
Theorem 4.6 and Remark 4.7 in \cite{JN94}.


\end{document}